\documentclass[11pt,a4paper]{article}

\usepackage{amssymb}
\usepackage{amsmath}
\usepackage{epic}
\usepackage{amsthm}
\usepackage{amscd}
\usepackage[all,2cell]{xy}
\usepackage[dvips]{graphics}

\newtheorem{Thm}{Theorem}[section]
\newtheorem{Cor}[Thm]{Corollary}
\newtheorem{Lem}[Thm]{Lemma}
\newtheorem{Prop}[Thm]{Proposition}

\theoremstyle{remark}

\newtheorem*{pf}{\rm\textbf{Proof}}
\newtheorem{eg}{\rm\textbf{Example}}

\begin{document}

\title{\Large \bf On the Affine Schur Algebra of Type $A$\ \ \ II }
\date{}
\author{\large Dong Yang\footnote{The author acknowledges partial support by the AsiaLink network
Algebras and Representations in China and Europe,
ASI/B7-301/98/679-11, and by National Natural Science Foundation
of China No.10131010.}}

\maketitle

\abstract{\begin{quote} By studying certain kind of centralizer
algebras of the affine Schur algebra $\widetilde{S}(n,r)$ we show
that
$\widetilde{S}(n,r)$ is Noetherian and we determine its center. Assuming $n\geq r$, we show that $\widetilde{S}(n+1,r)$ is Morita equivalent to $\widetilde{S}(n,r)$, and the Schur functor is an equivalence under certain conditions.
\end{quote}}

\section{Introduction}\label{S:introduction}

The affine ($q$-)Schur algebra (of type $A$) has been studied
by~\cite{GV}~\cite{RMG}~\cite{L}~\cite{L2}~\cite{Y}~\cite{VV},
which provide various equivalent definitions of the algebra.

In this paper we define the affine Schur algebra
$\widetilde{S}(n,r)$ ($n,r\in\mathbb{N}$) by giving a basis and
the structure constants. We investigate certain centralizer
algebras of $\widetilde{S}(n,r)$ of the form
$e\widetilde{S}(n,r)e$, where $e$ is an idempotent of
$\widetilde{S}(n,r)$. The most interesting ones are those
$\widetilde{S}(\underline{i})\xi_{\underline{i},\underline{i}}\widetilde{S}(n,r)\xi_{\underline{i},\underline{i}}$
 for $\underline{i}\in
I(n,r)$, and those isomorphic to $\widetilde{S}(n',r)$ for some
$n'\in\mathbb{N}$. As applications, we obtain the following
results. Firstly, when $n\geq r$, $\widetilde{S}(n+1,r)$ is Morita
equivalent to $\widetilde{S}(n,r)$ (Theorem~\ref{T:minus1}).
Secondly, when $n\geq r$, the Schur functor is well-defined and it
is an equivalence --- thus the affine Schur algebra is Morita
equivalent to the group algebra of the extended affine Weyl group
--- provided the characteristic of the base field is $0$ or greater
than $r$ (Theorem~\ref{T:schurfunctor}). These two equivalences
are affine analogues to the results in the finite case ({\em
cf.}~\cite{JAG}). Thirdly, $\widetilde{S}(n,r)$ is Noetherian
(Theorem~\ref{T:noetherian}). Besides, we determine the center of
$\widetilde{S}(n,r)$. More precisely, its center is isomorphic to
the algebra $K[t_{1},\cdots,t_{r-1},t_{r},t_{r}^{-1}]$ where
$t_{1},\cdots,t_{r}$ are indeterminates (Theorem~\ref{T:center}).

This paper is organized as follows. In
Section~\ref{S:affineschuralgebra}, we give the definition of the
affine Schur algebra and recall some basic properties.
Section~\ref{S:certaincentralizeralgebras} is devoted to the study
of certain centralizer algebras. In this section we show the
Morita equivalences stated above and prove that the affine Schur
algebra is Noetherian. In Section~\ref{S:center} we determine the
center. Section~\ref{S:examples} provides some examples of the
affine Schur algebra.

\section{The affine Schur algebra}\label{S:affineschuralgebra}
First let us introduce the setting.

$K$ will be an infinite field, and $n,r\in\mathbb{N}$. Let
$\Sigma_{r}$ denote the symmetric group on $r$ letters and
$\widehat{\Sigma}_{r}=\Sigma_{r}\ltimes \mathbb{Z}^{r}$ the
extended affine Weyl group of type $A_{r-1}$.

For a set $S$, we denote by $I(S,r)$ the set
$\{\underline{i}=(i_{1},\cdots,i_{r})\ |\ i_{t}\in S,\
t=1,\cdots,r\}$ of all $r$-tuples of elements in $S$. We often
omit the brackets and commas in the expression of
$\underline{i}\in I(S,r)$ if it does not cause confusion. Then
$\Sigma_{r}$ acts on the right on $I(S,r)$ by place permutation.
We will abbreviate $I(\{1,\cdots,n\},r)$ to $I(n,r)$. $\Sigma_{n}$
acts on the left on $I(n,r)$.

$\widehat{\Sigma}_{r}$ acts on $I(\mathbb{Z},r)$ on the right with
$\Sigma_{r}$ acting by place permutation and $\mathbb{Z}^{r}$
acting by shifting, i.e. $\underline{i}\varepsilon=\underline{i}+
n \varepsilon$ for $\underline{i}\in I(\mathbb{Z},r)$ and
$\varepsilon\in\mathbb{Z}^{r}$, and on $I(\mathbb{Z},r)\times
I(\mathbb{Z},r)$ diagonally. This action depends on the number
$n$. Note that a representative set for $I(n,r)/\Sigma_{r}$ is
also a representative set for
$I(\mathbb{Z},r)/\widehat{\Sigma}_{r}$. For $\underline{i}\in
I(\mathbb{Z},r)$, set $\widehat{\Sigma}_{\underline{i}}$ (resp.
$\Sigma_{\underline{i}}$) be the stabilizer group of
$\underline{i}$ in $\widehat{\Sigma}_{r}$ (resp. $\Sigma_{r}$),
and for $\underline{i},\underline{j}\in I(\mathbb{Z},r)$, set
$\widehat{\Sigma}_{\underline{i},\underline{j}}=\widehat{\Sigma}_{\underline{i}}\cap
\widehat{\Sigma}_{\underline{j}}$ (resp.
$\Sigma_{\underline{i},\underline{j}}=\Sigma_{\underline{i}}\cap\Sigma_{\underline{j}}$),
and so on. Note that if $\underline{i}\in I(n,r)$ then
$\widehat{\Sigma}_{\underline{i}}=\Sigma_{\underline{i}}$.

\vskip5pt

To each pair $(\underline{i},\underline{j})\in
I(\mathbb{Z},r)\times I(\mathbb{Z},r)$, we associate an element
$\xi_{\underline{i},\underline{j}}$ such that
$\xi_{\underline{i},\underline{j}}=\xi_{\underline{k},\underline{l}}$
if and only if
$(\underline{i},\underline{j})\sim_{\widehat{\Sigma}_{r}}(\underline{k},\underline{l})$.
{\em The affine Schur algebra $\widetilde{S}(n,r)$} is defined to
be the $K$-algebra with basis
$\{\xi_{\underline{i},\underline{j}}|\underline{i},\underline{j}\in
I(\mathbb{Z},r)\}$ and multiplication given by
\[
\xi_{\underline{i},\underline{j}}\xi_{\underline{k},\underline{l}}=\sum_{(\underline{p},\underline{q})\in
I(\mathbb{Z},r)\times
I(\mathbb{Z},r)/\widehat{\Sigma}_{r}}{Z(\underline{i},\underline{j},\underline{k},\underline{l},\underline{p},\underline{q})\xi_{\underline{p},\underline{q}}}
\]
where
$Z(\underline{i},\underline{j},\underline{k},\underline{l},\underline{p},\underline{q})=
\#\{\underline{s}\in
I(\mathbb{Z},r)|(\underline{i},\underline{j})\sim_{\widehat{\Sigma}_{r}}(\underline{p},\underline{s}),(\underline{s},\underline{q})\sim_{\widehat{\Sigma}{r}}(\underline{k},\underline{l})\}$.
We have

\begin{Prop}\label{P:basic} (~\cite{Y}Proposition4.2)

{\rm(i)}
 $\xi_{\underline{i},\underline{j}}\xi_{\underline{k},\underline{l}}=0$
unless $\underline{j} \sim_{\widehat{\Sigma}_{r}} \underline{k}$.

{\rm(ii)}
$\xi_{\underline{i},\underline{i}}\xi_{\underline{i},\underline{j}}=\xi_{\underline{i},\underline{j}}=\xi_{\underline{i},\underline{j}}\xi_{\underline{j},\underline{j}}$,
for $\underline{i},\underline{j}\in I(\mathbb{Z},r)$.

{\rm(iii)} $\sum_{\underline{i}\in
I(n,r)/\Sigma_{r}}{\xi_{\underline{i},\underline{i}}}$ is a
decomposition of the identity into orthogonal idempotents.

\end{Prop}

We have another product formula, which is proved in the end of
~\cite{Y} Section4.

\begin{Prop}\label{P:productformula} For $\underline{i},\underline{j},\underline{l}\in
I(\mathbb{Z},r)$, we have
\[\xi_{\underline{i},\underline{j}}\xi_{\underline{j},\underline{l}}=\sum_{\delta\in
\widehat{\Sigma}_{\underline{j},\underline{l}}\backslash\widehat{\Sigma}_{\underline{j}}/\widehat{\Sigma}_{\underline{i},\underline{j}}}[\widehat{\Sigma}_{\underline{i},\underline{l}\delta}:\widehat{\Sigma}_{\underline{i},\underline{j},\underline{l}
\delta}]\xi_{\underline{i},\underline{l}\delta}=\sum_{\delta\in
\widehat{\Sigma}_{\underline{i},\underline{j}}\backslash\widehat{\Sigma}_{\underline{j}}/\widehat{\Sigma}_{\underline{j},\underline{l}}}[\widehat{\Sigma}_{\underline{i}\delta,\underline{l}}:\widehat{\Sigma}_{\underline{i}\delta,\underline{j},\underline{l}}]\xi_{\underline{i}\delta,\underline{l}}\
.\]

\end{Prop}

Let \ $\bar{}: \mathbb{Z}\rightarrow \{1,\cdots,n\}$ be the map
taking least positive remainder modulo $n$. It can be extended to
\ $\bar{}: I(\mathbb{Z},r)\rightarrow I(n,r)$. Note that
$\xi_{\underline{i},\underline{j}}=\xi_{\overline{\underline{i}},\underline{j}+n\varepsilon^{1}}=\xi_{\underline{i}+n\varepsilon^{2},\overline{\underline{j}}}$\
, where
$\varepsilon^{1}=\frac{\overline{\underline{i}}-\underline{i}}{n}$,
$\varepsilon^{2}=\frac{\overline{\underline{j}}-\underline{j}}{n}$
are both in $\mathbb{Z}^{r}$. Since
$\{(\underline{i},\underline{j}+n\varepsilon)\ |\ \underline{i}\in
I(n,r)/\Sigma_{r}, \underline{j}\in I(n,r)/\Sigma_{\underline{i}},
\varepsilon\in
\mathbb{Z}^{r}/\Sigma_{\underline{i},\underline{j}}\}$ is a set of
representatives of the $\widehat{\Sigma}_{r}$-orbits of
$I(\mathbb{Z},r)\times I(\mathbb{Z},r)$, the set
$\{\xi_{\underline{i},\underline{j}}\ |\
\underline{i},\underline{j}\in
I(\mathbb{Z},r)\}=\{\xi_{\underline{i},\underline{j}}\ |\
(\underline{i},\underline{j})\in I(\mathbb{Z},r)\times
I(\mathbb{Z},r)/\widehat{\Sigma}_{r}\}$ equals the set
$\{\xi_{\underline{i},\underline{j}+n\varepsilon}\ |\
\underline{i},\underline{j}\in I(n,r), \varepsilon\in
\mathbb{Z}^{r}\}=\{\xi_{\underline{i},\underline{j}+n\varepsilon}\
|\ \underline{i}\in I(n,r)/\Sigma_{r}, \underline{j}\in
I(n,r)/\Sigma_{\underline{i}}, \varepsilon\in
\mathbb{Z}^{r}/\Sigma_{\underline{i},\underline{j}}\}$. Thus we
can rewrite Proposition~\ref{P:productformula} as follows.

\begin{Prop}\label{P:productformula2} For $\underline{i},\underline{j},\underline{l}\in
I(n,r)$, $\varepsilon,\varepsilon'\in\mathbb{Z}^{r}$, we have
\[\begin{array}{c}\xi_{\underline{i},\underline{j}+n\varepsilon}\xi_{\underline{j},\underline{l}+n\varepsilon'}
=\sum_{\delta\in
{\Sigma}_{\underline{j},\underline{l},\varepsilon'}\backslash
{\Sigma}_{\underline{j}}/{\Sigma}_{\underline{i},\underline{j},\varepsilon}}[{\Sigma}_{\underline{i},\underline{l}\delta,\varepsilon'\delta+\varepsilon}:{\Sigma}_{\underline{i},\underline{j},\underline{l}\delta,\varepsilon'\delta,\varepsilon}]\xi_{\underline{i},\underline{l}\delta+n(\varepsilon'\delta+\varepsilon)}\\
=\sum_{\delta\in
{\Sigma}_{\underline{i},\underline{j},\varepsilon}\backslash
{\Sigma}_{\underline{j}}/{\Sigma}_{\underline{j},\underline{l},\varepsilon'}}[{\Sigma}_{\underline{i}\delta,\underline{l},\varepsilon'+\varepsilon\delta}:{\Sigma}_{\underline{i}\delta,\underline{j},\underline{l},\varepsilon',\varepsilon\delta}]\xi_{\underline{i}\delta,\underline{l}+n(\varepsilon'+\varepsilon\delta)}\
.\end{array}\]
\end{Prop}
\vskip20pt

\section{Certain centralizer
algebras}\label{S:certaincentralizeralgebras} In this section we
study centralizer algebras of the form $e\widetilde{S}(n,r)e$,
where $e$ is an idempotent of $\widetilde{S}(n,r)$. We show a few
Morita equivalences and prove that $\widetilde{S}(n,r)$ is
Noetherian.

Proposition~\ref{P:basic}(iii) says that $\sum_{\underline{i}\in
I(n,r)/\Sigma_{r}}{\xi_{\underline{i},\underline{i}}}$ is a
decomposition of the identity into orthogonal idempotents.
Correspondingly $\oplus_{\underline{i}\in
I(n,r)/\Sigma_{r}}\widetilde{S}(n,r)\xi_{\underline{i},\underline{i}}$
is a decomposition of the regular $\widetilde{S}(n,r)$-module into
a direct sum of projective modules.

\begin{Lem}\label{L:summand}
Let $\underline{i},\underline{j}\in I(n,r)$ with
$\Sigma_{\underline{i}}\geq \Sigma_{\underline{j}}$. Assume
$charK\dagger [\Sigma_{\underline{i}}:\Sigma_{\underline{j}}]$,
then as an $\widetilde{S}(n,r)$-module
$\widetilde{S}(n,r)\xi_{\underline{i},\underline{i}}$ is a direct
summand of $\widetilde{S}(n,r)\xi_{\underline{j},\underline{j}}$.
\end{Lem}
\begin{pf}
Note that under that assumption
\[
\xi_{\underline{i},\underline{j}}\xi_{\underline{j},\underline{i}}=\sum_{\delta\in\Sigma_{\underline{i},\underline{j}}\backslash\Sigma_{\underline{j}}/\Sigma_{\underline{i},\underline{j}}}[\Sigma_{\underline{i},\underline{i}\delta}:\Sigma_{\underline{i},\underline{i}\delta,\underline{j}}]\xi_{\underline{i},\underline{i}\delta}=[\Sigma_{\underline{i}}:\Sigma_{\underline{j}}]\xi_{\underline{i},\underline{i}}
\] So
$\phi:\widetilde{S}(n,r)\xi_{\underline{i},\underline{i}}\rightarrow\widetilde{S}(n,r)\xi_{\underline{j},\underline{j}}$,
$\xi\mapsto \xi\xi_{\underline{i},\underline{j}}$, and
$\psi:\widetilde{S}(n,r)\xi_{\underline{j},\underline{j}}\rightarrow\widetilde{S}(n,r)\xi_{\underline{i},\underline{i}}$,
$\xi\mapsto\frac{\xi\xi_{\underline{j},\underline{i}}}{[\Sigma_{\underline{i}}:\Sigma_{\underline{j}}]}$
are both homomorphisms of $\widetilde{S}(n,r)$-modules and
$\psi\circ\phi=id$. Consequently we get the desired result.
$\square$
\end{pf}

For $\underline{i}\in I(n,r)$, denote by
$\widetilde{S}(\underline{i})$ the centralizer algebra
$\xi_{\underline{i},\underline{i}}\widetilde{S}(n,r)\xi_{\underline{i},\underline{i}}$.
The following is a corollary of Lemma~\ref{L:summand}.

\begin{Lem}\label{L:centrasubalg}
Let $\underline{i}\in I(n,r)$ and
$\underline{j}=\sigma(\underline{i})$ for some
$\sigma\in\Sigma_{n}$. Then
$\widetilde{S}(n,r)\xi_{\underline{i},\underline{i}}$ and
$\widetilde{S}(n,r)\xi_{\underline{j},\underline{j}}$ are
isomorphic $\widetilde{S}(n,r)$-modules. Consequently,
$\widetilde{S}(\underline{i})$ and $\widetilde{S}(\underline{j})$
are isomorphic $K$-algebras.
\end{Lem}
\begin{pf} The first statement follows from Lemma~\ref{L:summand}
since $\Sigma_{\underline{j}}=\Sigma_{\underline{i}}$. The second
statement is because isomorphic modules have isomorphic
endomorphism algebras. $\square$
\end{pf}

The next proposition shows that some centralizer algebras are in
fact affine Schur algebras with different parameters.

\begin{Prop}\label{P:centrasubalg} Let $N$ be a subset of $\{1,\cdots,n\}$, and
$\xi_{N}=\sum_{\underline{i}\in
I(N,r)/\Sigma_{r}}\xi_{\underline{i},\underline{i}}$. Then
$\xi_{N}\widetilde{S}(n,r)\xi_{N}\cong \widetilde{S}(\#N,r)$.
\end{Prop}
\begin{pf} Let $s=\# N$, and $N^{0}=\{1,\cdots,s\}$. Then there
exists a $\sigma\in\Sigma_{n}$ such that $\sigma$ induces a
bijection from $I(N,r)$ to $I(N^{0},r)$ given by
$\underline{i}\mapsto \sigma(\underline{i})$. Therefore by
Lemma~\ref{L:centrasubalg} $\widetilde{S}(n,r)\xi_{N}$ and
$\widetilde{S}(n,r)\xi_{N^{0}}$ are isomorphic
$\widetilde{S}(n,r)$-modules. As a result, their endomorphism
algebras $\xi_{N}\widetilde{S}(n,r)\xi_{N}$ and
$\xi_{N^{0}}\widetilde{S}(n,r)\xi_{N^{0}}$ are isomorphic
$K$-algebras. Finally, sending
$\xi_{\underline{i},\underline{j}+n\varepsilon}$ to
$\xi_{\underline{i},\underline{j}+s\varepsilon}$ is a $K$-algebra
isomorphism from $\xi_{N^{0}}\widetilde{S}(n,r)\xi_{N^{0}}$ to
$\widetilde{S}(s,r)$. $\square$

\end{pf}

Now we are able to establish a Morita equivalence.

\begin{Thm}\label{T:minus1} Assume $n\geq r$, then
$\widetilde{S}(n+1,r)$ is Morita equivalent to
$\widetilde{S}(n,r)$.
\end{Thm}
\begin{pf} Since $n\geq r$, we have that for any $\underline{i}\in I(n+1,r)$
there exists $\sigma\in \Sigma_{n+1}$ such that
$\sigma(\underline{i})\in I(n,r)$. Therefore
$\widetilde{S}(n+1,r)\xi_{N}$ is a progenerator where
$N=\{1,\cdots,n\}$. In particular, $\widetilde{S}(n+1,r)$ is
Morita equivalent to $\xi_{N}\widetilde{S}(n+1,r)\xi_{N}$, which
is isomorphic to $\widetilde{S}(n,r)$ by
Proposition~\ref{P:centrasubalg}. $\square$
\end{pf}

Next we will concentrate on the algebras
$\widetilde{S}(\underline{i})$. We start with a special case.

\subsection{Special case : $\widetilde{S}(1\cdots
1)$}\label{S:1...1}

$\widetilde{S}(1\cdots 1)=K\{\xi_{1\cdots 1,1\cdots
1+n\varepsilon}\ |\ \varepsilon\in\mathbb{Z}^{r}\}$. Following
Proposition~\ref{P:productformula2} we write down the product
formula for this algebra \[\begin{array}{c}\xi_{1\cdots 1,1\cdots
1+n\varepsilon}\xi_{1\cdots 1,1\cdots
1+n\varepsilon'}=\sum_{\delta\in\Sigma_{\varepsilon'}\backslash\sigma_{r}/\Sigma_{\varepsilon}}[\Sigma_{\varepsilon'\delta+\varepsilon}:\Sigma_{\varepsilon'\delta,\varepsilon}]\xi_{1\cdots
1,1\cdots
1+n(\varepsilon'\delta+\varepsilon)}\\
=\sum_{\delta\in\Sigma_{\varepsilon}\backslash\sigma_{r}/\Sigma_{\varepsilon'}}[\Sigma_{\varepsilon\delta+\varepsilon'}:\Sigma_{\varepsilon\delta,\varepsilon'}]\xi_{1\cdots
1,1\cdots 1+n(\varepsilon\delta+\varepsilon')}\
,\end{array}\]where $\varepsilon,\varepsilon'\in\mathbb{Z}^{r}$.
In particular $\widetilde{S}(1\cdots 1)$ is a commutative algebra.
Since $\Sigma_{1\cdots 1}=\Sigma_{r}$, we see that $\xi_{1\cdots
1,1\cdots 1+n\varepsilon}=\xi_{1\cdots 1,1\cdots
1+n\varepsilon\sigma}$ for any $\varepsilon\in\mathbb{Z}^{r}$ and
$\sigma\in\Sigma_{r}$. Thus for a basis element $\xi_{1\cdots
1,1\cdots 1+n\varepsilon}$ we can always choose $\varepsilon$ to
be weakly decreasing. In other words, $\xi_{1\cdots 1,1\cdots
1+n\varepsilon}$, $\varepsilon\in\mathbb{Z}^{r}$ weakly
decreasing, form a basis for $\widetilde{S}(1\cdots 1)$.

We denote by $\widetilde{S}(1\cdots 1)^{+}$ the subalgebra of
$\widetilde{S}(1\cdots 1)$ with basis $\{\xi_{1\cdots 1,1\cdots
1+n\varepsilon}\ |\
\varepsilon=(\varepsilon_{1},\cdots,\varepsilon_{r})\in\mathbb{Z}^{r},\varepsilon_{1}\geq\cdots\geq\varepsilon_{r}\geq
0\}$. Then $\widetilde{S}(1\cdots 1)^{+}$ is a positively graded
$K$-algebra with grading given by $deg(\xi_{1\cdots 1,1\cdots
1+n(\varepsilon_{1},\cdots,\varepsilon_{r})})=\varepsilon_{1}+\cdots+\varepsilon_{r}$.
Let $e^{1},\cdots,e^{r}$ be the natural $\mathbb{Z}$-basis for
$\mathbb{Z}^{r}$. Let $\varepsilon^{k}=\sum_{i=1}^{k}e^{i}$ and
$t_{k}=\xi_{1\cdots 1,1\cdots 1+n\varepsilon^{k}}$ for
$k=1,\cdots,r$.

\begin{Prop}\label{P:1...1+} $\widetilde{S}(1\cdots 1)^{+}=K\xi_{1\cdots 1,1\cdots 1}[t_{1},\cdots,t_{r-1},t_{r}]$ is a polynomial algebra in $r$ indeterminates.
\end{Prop}
\begin{pf}
We prove in two steps: $t_{1},\cdots,t_{r}$ generate
$\widetilde{S}(1\cdots 1)^{+}$ and they are algebraically
independent. For the first step, it suffices to show basis
elements $\xi_{1\cdots 1,1\cdots 1+n\varepsilon}$,
$\varepsilon\in\mathbb{Z}^{r}$ weakly decreasing, are generated by
them. Induct on the degree and the number of $s$'s with
$\varepsilon_{s}=0$.

If $\varepsilon_{r}\not=0$, then $\xi_{1\cdots 1,1\cdots
1+n\varepsilon} =\xi_{1\cdots 1,1\cdots
1+n(\varepsilon-\varepsilon_{r}\varepsilon^{r})}t_{r}^{\varepsilon_{r}}$
is the product of an element of less degree and a power of
$t_{r}$.

Assume $\varepsilon_{k+1}=\cdots=\varepsilon_{r}=0$ but
$\varepsilon_{k}\not= 0$, then
\[
\xi_{1\cdots 1,1\cdots 1+n\varepsilon}=\xi_{1\cdots 1,1\cdots
1+n(\varepsilon-\varepsilon^{k})}t_{k}
-\sum_{\delta}[\Sigma_{\varepsilon^{k}\delta+\varepsilon-\varepsilon^{k}}:\Sigma_{\varepsilon^{k}\delta,\varepsilon-\varepsilon^{k}}]\xi_{1\cdots
1,1\cdots 1+n(\varepsilon^{k}\delta+\varepsilon-\varepsilon^{k})}
\]
where the sum is over a set of representatives of all non-trivial
double cosets
$\Sigma_{\varepsilon^{k}}\delta\Sigma_{\varepsilon-\varepsilon^{k}}$
of $\Sigma_{r}$. Note that $\xi_{1\cdots 1,1\cdots
1+n(\varepsilon-\varepsilon^{k})}$ is of less degree and
$\varepsilon^{k}\delta+\varepsilon-\varepsilon^{k}$ has less zero
entries than $\varepsilon$. By induction we finish this step.

Now let us prove the second step. The $K$-dimension of the
homogeneous component of $\widetilde{S}(1\cdots 1)^{+}$ of degree
$nm$ is
$\#\{\varepsilon\in\mathbb{Z}^{r}|\varepsilon_1\geq\cdots\geq\varepsilon_{r}\geq
0, \varepsilon_{1}+\cdots+\varepsilon_{r}=m\}$. This set is in
bijection with the set
$\{(n_{1},\cdots,n_{r})\in\mathbb{Z}^{r}|n_{1},\cdots,n_{r}\geq
0,n_{1}+2n_{2}+\cdots+rn_{r}=m\}$, whose cardinality equals the
$K$-dimension of the homogeneous component of the polynomial
algebra $K[x_{1},\cdots,x_{r}]$ of degree $nm$ with
$deg(x_{k})=nk$, $k=1,\cdots,r$. This dimension comparison shows
that $t_{1},\cdots,t_{r}$ are algebraically independent.
 $\square$
\end{pf}

Let $\varepsilon\in\mathbb{Z}^{r}$ be weakly decreasing. If
$\varepsilon_{r}<0$, then $ \xi_{1\cdots 1,1\cdots 1+n\varepsilon}
=\xi_{1\cdots 1,1\cdots
1+n(\varepsilon-\varepsilon_{r}\varepsilon^{r})}(t_{r}^{-1})^{-\varepsilon_{r}}
$ is the product of an element $\xi_{1\cdots 1,1\cdots
1+n(\varepsilon-\varepsilon_{r}\varepsilon^{r})}$ in
$\widetilde{S}(1\cdots 1)^{+}$ and a power of $t_{r}^{-1}$, where
$t_{r}^{-1}=\xi_{1\cdots 1,1\cdots 1+n(-\varepsilon^{r})}$.

\begin{Prop}\label{P:1...1}

{\rm (i)} As a $K$-algebra $\widetilde{S}(1\cdots 1)$ is
isomorphic to $ K[t_{1},\cdots,t_{r-1},t_{r},t_{r}^{-1}]$, where
$t_{1},\cdots,t_{r}$ are indeterminates. In particular, the affine
Schur algebra  $\widetilde{S}(1,r)$ is isomorphic to $
K[t_{1},\cdots,t_{r-1},t_{r},t_{r}^{-1}]$.

{\rm (ii)} $\xi_{1\cdots 1,1\cdots 1}$ is a primitive idempotent
of $\widetilde{S}(n,r)$.
\end{Prop}

\subsection{A commutative subalgebra of  $\widetilde{S}(\underline{i})$}

Fix $\underline{i}\in I(n,r)$ in this subsection. Then
$\widetilde{S}(\underline{i})=K\{\xi_{\underline{i},\underline{i}w}\
|\
w\in\widehat{\Sigma}_{r}\}=K\{\xi_{\underline{i},\underline{i}\sigma+n\varepsilon}\
|\ \sigma\in\Sigma_{r},\varepsilon\in\mathbb{Z}^{r}\}$. The rest
of this subsection is devoted to the study of a commutative
subalgebra of $\widetilde{S}(\underline{i})$.

Let
$B_{\underline{i}}=K\{\xi_{\underline{i},\underline{i}+n\varepsilon}\}$.
By Proposition~\ref{P:productformula2} for
$\varepsilon,\varepsilon'\in\mathbb{Z}^{r}$,
\[\begin{array}{c}\xi_{\underline{i},\underline{i}+n\varepsilon}\xi_{\underline{i},\underline{i}+n\varepsilon'}=\sum_{\delta\in\Sigma_{\underline{i},\varepsilon}\backslash\Sigma_{\underline{i}}/\Sigma_{\underline{i},\varepsilon'}}[\Sigma_{\underline{i},\varepsilon'+\varepsilon\delta}:\Sigma_{\underline{i},\varepsilon',\varepsilon\delta}]\xi_{\underline{i},\underline{i}+n(\varepsilon'+\varepsilon\delta)}\\
=\sum_{\delta\in\Sigma_{\underline{i},\varepsilon'}\backslash\Sigma_{\underline{i}}/\Sigma_{\underline{i},\varepsilon}}[\Sigma_{\underline{i},\varepsilon+\varepsilon'\delta}:\Sigma_{\underline{i},\varepsilon,\varepsilon'\delta}]\xi_{\underline{i},\underline{i}+n(\varepsilon+\varepsilon'\delta)}\
. \end{array}\] It follows that $B_{\underline{i}}$ is a
commutative subalgebra of $\widetilde{S}(\underline{i})$.
Moreover,

\begin{Prop}\label{P:bimaximal} $B_{\underline{i}}$ is a maximal commutative subalgebra of
$\widetilde{S}(\underline{i})$.
\end{Prop}
\begin{pf} Let $\xi=\sum_{\sigma\in\Sigma_{\underline{i}}\backslash\Sigma_{r}/\Sigma_{\underline{i}}}\sum_{\varepsilon\in \mathbb{Z}^{r}/\Sigma_{\underline{i},\underline{i}\sigma}}\lambda_{\sigma,\varepsilon}\xi_{\underline{i},\underline{i}\sigma+n\varepsilon}\in \widetilde{S}(\underline{i})$  be such
that for any $\varepsilon^{0}\in\mathbb{Z}^{r}$,
\[\xi\xi_{\underline{i},\underline{i}+n\varepsilon^{0}}=\xi_{\underline{i},\underline{i}+n\varepsilon^{\tiny 0}}\xi\ ,\]
\[i.e.\ \ \sum_{\sigma\in\Sigma_{\underline{i}}\backslash\Sigma_{r}/\Sigma_{\underline{i}}}\ \sum_{\varepsilon\in
\mathbb{Z}^{r}/\Sigma_{\underline{i},\underline{i}\sigma}}\lambda_{\sigma,\varepsilon}\sum_{\delta\in\Sigma_{\underline{i}\sigma,\varepsilon^{0}\sigma}\backslash\Sigma_{\underline{i}\sigma}/\Sigma_{\underline{i}\sigma,\underline{i},\varepsilon}}[\Sigma_{\underline{i},\underline{i}\sigma,\varepsilon+\varepsilon^{0}\sigma\delta}:\Sigma_{\underline{i},\underline{i}\sigma,\varepsilon,\varepsilon^{0}\sigma\delta}]\xi_{\underline{i},\underline{i}\sigma+n(\varepsilon+\varepsilon^{0}\sigma\delta)}\]
\[
=\sum_{\sigma\in\Sigma_{\underline{i}}\backslash\Sigma_{r}/\Sigma_{\underline{i}}}\
\sum_{\varepsilon\in
\mathbb{Z}^{r}/\Sigma_{\underline{i},\underline{i}\sigma}}\lambda_{\sigma,\varepsilon}\sum_{\delta\in\Sigma_{\underline{i},\varepsilon^{0}}\backslash\Sigma_{\underline{i}}/\Sigma_{\underline{i},\underline{i}\sigma,\varepsilon}}[\Sigma_{\underline{i},\underline{i}\sigma,\varepsilon+\varepsilon^{0}\delta}:\Sigma_{\underline{i},\underline{i}\sigma,\varepsilon,\varepsilon^{0}\delta}]\xi_{\underline{i},\underline{i}\sigma+n(\varepsilon+\varepsilon^{0}\delta)}\
\ \ \ \ (*)\] \vskip5pt
 Take $t\in\mathbb{N}$ large enough, say,
$t=10\times max\{|\varepsilon_{s}|\ \big|\ s=1,\cdots,r,
\varepsilon : \lambda_{\sigma,\varepsilon}\neq 0 {\text\rm\ for\
some\ } \sigma\}$ and take $\varepsilon^{0}=t(1,2,\cdots,r)$. Then
for
${\sigma\in\Sigma_{\underline{i}}\backslash\Sigma_{r}/\Sigma_{\underline{i}}}$,
${\varepsilon\in
\mathbb{Z}^{r}/\Sigma_{\underline{i},\underline{i}\sigma}}$,
$\delta\in\Sigma_{\underline{i},\varepsilon^{0}}\backslash\Sigma_{\underline{i}}/\Sigma_{\underline{i},\underline{i}\sigma,\varepsilon}$,
we have $\Sigma_{\varepsilon+\varepsilon^{0}\delta}=1$. In
particular,
$[\Sigma_{\underline{i},\underline{i}\sigma,\varepsilon+\varepsilon^{0}\delta}:\Sigma_{\underline{i},\underline{i}\sigma,\varepsilon,\varepsilon^{0}\delta}]=1$.
Moreover, for
${\sigma'\in\Sigma_{\underline{i}}\backslash\Sigma_{r}/\Sigma_{\underline{i}}}$,
${\varepsilon'\in
\mathbb{Z}^{r}/\Sigma_{\underline{i},\underline{i}\sigma'}}$,
${\delta'\in\Sigma_{\underline{i},\varepsilon^{0}}\backslash\Sigma_{\underline{i}}/\Sigma_{\underline{i},\underline{i}\sigma',\varepsilon'}}$,
we have
$\xi_{\underline{i},\underline{i}\sigma+n(\varepsilon+\varepsilon^{0}\delta)}=\xi_{\underline{i},\underline{i}\sigma'+n(\varepsilon'+\varepsilon^{0}\delta')}$
implies
$(\underline{i},\underline{i}\sigma+n(\varepsilon+\varepsilon^{0}\delta))\sim_{\widehat{\Sigma}_{r}}(\underline{i},\underline{i}\sigma'+n(\varepsilon'+\varepsilon^{0}\delta'))$,
i.e. there exists $\tau\in\widehat{\Sigma}_{r}$ such that
$\underline{i}=\underline{i}\tau, \
\underline{i}\sigma=\underline{i}\sigma'\tau,\
\varepsilon=\varepsilon'\tau,\
\varepsilon^{0}\delta=\varepsilon^{0}\delta'\tau$. Thus
$\tau\in\Sigma_{\underline{i}}$,
$\sigma'\in\Sigma_{\underline{i}}\sigma\tau^{-1}$, and
$\delta'\tau=\delta$ since $\Sigma_{\varepsilon^{0}}$ is trivial.
Therefore $\sigma=\sigma'$, and hence
$\tau\in\Sigma_{\underline{i},\underline{i}\sigma}$. So
$\varepsilon=\varepsilon'$, and then
$\tau\in\Sigma_{\underline{i},\underline{i}\sigma,\varepsilon}$,
and hence $\delta=\delta'$.

Now suppose
${\sigma\in\Sigma_{\underline{i}}\backslash\Sigma_{r}/\Sigma_{\underline{i}}},{\varepsilon\in
\mathbb{Z}^{r}/\Sigma_{\underline{i},\underline{i}\sigma}},{\delta\in\Sigma_{\underline{i},\varepsilon^{0}}\backslash\Sigma_{\underline{i}}/\Sigma_{\underline{i},\underline{i}\sigma,\varepsilon}}$
satisfy $\lambda_{\sigma,\varepsilon}\neq 0$. Then by $(*)$ there
exist
${\sigma'\in\Sigma_{\underline{i}}\backslash\Sigma_{r}/\Sigma_{\underline{i}}},{\varepsilon'\in
\mathbb{Z}^{r}/\Sigma_{\underline{i},\underline{i}\sigma'}},{\delta'\in\Sigma_{\underline{i}\sigma,\varepsilon^{0}\sigma}\backslash\Sigma_{\underline{i}\sigma}/\Sigma_{\underline{i},\underline{i}\sigma,\varepsilon}}$
such that
$\xi_{\underline{i},\underline{i}\sigma'+n(\varepsilon'+\varepsilon^{0}\sigma'\delta')}=\xi_{\underline{i},\underline{i}\sigma+n(\varepsilon+\varepsilon^{0}\delta)}$,
 i.e. there exists $\tau\in\widehat{\Sigma}_{r}$ such that
$\underline{i}=\underline{i}\tau,\
\underline{i}\sigma'=\underline{i}\sigma\tau,\
\varepsilon'=\varepsilon\tau,\
\varepsilon^{0}\sigma'\delta'=\varepsilon^{0}\delta\tau$. Similar
arguments as above show that $\tau\in
{\Sigma}_{\underline{i},\underline{i}\sigma,\varepsilon}$, and
$\sigma'=\sigma$, $\varepsilon=\varepsilon'$,
$\sigma\delta'=\delta\tau$. This implies
$\sigma=(\sigma\delta'\sigma^{-1})^{-1}\delta\tau\in\Sigma_{\underline{i}}$,
i.e. $\sigma$ is trivial. In a word,
$\lambda_{\sigma,\varepsilon}\neq 0$ implies $\sigma$ is trivial.
That is, $\xi\in B_{\underline{i}}$.
 $\square$
\end{pf}

 Let
$\lambda=(\lambda_{1},\cdots,\lambda_{n})$ be the weight of
$\underline{i}$, i.e. $\lambda_{s}=\#\{k\ |\ k=1,\cdots,r,\
i_{k}=s\}$, and let $\alpha$ be the number of nonzero entries of
$\lambda$. Then

\begin{Prop}\label{P:commsubalg}
$B_{\underline{i}}$ is a polynomial algebra in $r-\alpha$
indeterminates over a Laurent polynomial algebra in $\alpha$
indeterminates over $K$. In particular, $B_{\underline{i}}$ is
Noetherian.
\end{Prop}
\begin{pf}
For $s=1,\cdots,n$, set $\underline{i}(s,r)=s\cdots s\in I(n,r)$.
We may assume
$\underline{i}=\underline{i}(1,\lambda_{1})\cdots\underline{i}(n,\lambda_{n})$.
For $s=1,\cdots,n$, and $\varepsilon\in\mathbb{Z}^{r}$, set
\[\begin{array}{c}\Phi_{s}=\{0\}^{\lambda_{1}+\cdots+\lambda_{s-1}}\times
\mathbb{Z}^{\lambda_{s}}\times
\{0\}^{\lambda_{s+1}+\cdots+\lambda_{n}},\\
\theta^{0}_{s}(\varepsilon)=(\varepsilon_{\lambda_{1}+\cdots+\lambda_{s-1}+1},\cdots,\varepsilon_{\lambda_{1}+\cdots+\lambda_{s}})\in\mathbb{Z}^{\lambda_{s}},\\
\theta_{s}(\varepsilon)=(0,\cdots,0,\theta^{0}_{s}(\varepsilon),0,\cdots,0)\in
\Phi_{s}\ .\end{array}\]

Let $B_{s}=K\{\xi_{\underline{i},\underline{i}+n\varepsilon}\ |\
\varepsilon\in\Phi_{s}\}$. If $\lambda_{s}=0$ then $B_{s}\cong K$.
Otherwise
$\xi_{\underline{i},\underline{i}+n\varepsilon}\mapsto\xi_{1\cdots
1,1\cdots 1+\theta^{0}_{s}(\varepsilon)}$ defines an algebra
isomorphism $B_{s}\cong \widetilde{S}(1,\lambda_{s})$. By
Proposition~\ref{P:1...1} the latter algebra is a polynomial
algebra with $\lambda_{s}-1$ indeterminates over a Laurent
polynomial algebra in $1$ indeterminates. Now for
$\varepsilon\in\mathbb{Z}^{r}$, we have
$\xi_{\underline{i},\underline{i}+n\varepsilon}=\prod_{s=1}^{n}\xi_{\underline{i},\underline{i}+n\theta_{s}(\varepsilon)}$.
Therefore $B_{\underline{i}}\cong B_{1}\otimes\cdots\otimes
B_{n}$, and we are done. $\square$
\end{pf}

\subsection{$\widetilde{S}(n,r)$ is Noetherian}
In this subsection we will prove that $\widetilde{S}(n,r)$ is
Noetherian.

For $\varepsilon\in\mathbb{Z}^{r}$, let $\check{\varepsilon}$ be
the unique element in $\{\varepsilon\sigma\ |\
\sigma\in\Sigma_{r}\}$ such that
$\check{\varepsilon}_{1}\geq\check{\varepsilon}_{2}\geq\cdots\geq\check{\varepsilon}_{r}$.
For $\varepsilon,\varepsilon'\in\mathbb{Z}^{r}$, define
$\varepsilon> \varepsilon'$ if
$\varepsilon_{1}+\cdots+\varepsilon_{r}=\varepsilon'_{1}+\cdots+\varepsilon'_{r}$
and $\check{\varepsilon}>\check{\varepsilon'}$ according to the
lexicographic order. We way $\varepsilon\in\mathbb{Z}^{r}$ is {\em
successive} if $\{\varepsilon_{1},\cdots,\varepsilon_{r}\}$ is a
set of successive integers, $\varepsilon$ is {\em absolutely
successive} if in addition each entry is nonnegative and at least
one of them equals $0$.

Fix $\underline{i},\underline{j}\in I(n,r)$. By
Proposition~\ref{P:productformula2} the $K$-space $M$ spanned by
$\{\xi_{\underline{i},\underline{j}+n\varepsilon}\ |\
\varepsilon\in\mathbb{Z}^{r}\}$ is a $B_{\underline{i}}$-module.

\begin{Prop}\label{P:fg} The module $M$ is generated over $B_{\underline{i}}$ by
\[\{\xi_{\underline{i},\underline{j}+n\varepsilon}\ |\ \varepsilon {\text\rm\ is\ absolutely\ successive}\}\ .\]
In particular, it is finitely generated.
\end{Prop}
\begin{pf} We prove that the basis elements $\xi_{\underline{i},\underline{j}+n\varepsilon}$ is generated by the desired set.
Since
\[\xi_{\underline{i},\underline{j}+n\varepsilon}=\xi_{\underline{i},\underline{i}+nt(1,\cdots,1)}\xi_{\underline{i},\underline{j}+n(\varepsilon-t(1,\cdots,1))}\]
for $t=min\{\varepsilon_{1},\cdots,\varepsilon_{r}\}$, we may
assume $\varepsilon\geq 0$ and one of its entry is $0$. Induct on
the volume $\varepsilon_{1}+\cdots+\varepsilon_{r}$ of
$\varepsilon$.

If $\varepsilon$ satisfies
$\varepsilon_{1}+\cdots+\varepsilon_{r}=0$ or $1$, then
$\varepsilon$ is absolutely successive. Suppose
$\varepsilon_{1}+\cdots+\varepsilon_{r}\geq 2$. If $\varepsilon$
is successive, then it is absolutely successive. If $\varepsilon$
is not successive, then there is a partition
$\Omega\cup\Omega^{c}$ of $\{1,\cdots,r\}$ such that
$\varepsilon_{k}\geq\varepsilon_{k'}+2, \forall
k\in\Omega,k'\in\Omega^{c}$. Recall that $e^{1},\cdots,e^{r}$ is
the natural $\mathbb{Z}$-basis for $\mathbb{Z}^{r}$. Set
$e({\Omega})=\sum_{k\in\Omega}e^{k}$. Then
\[\xi_{\underline{i},\underline{j}+n\varepsilon}=\xi_{\underline{i},\underline{i}+ne_{\Omega}}\xi_{\underline{i},\underline{j}+n(\varepsilon-e(\Omega))}-\sum_{\delta}\mu_{\delta}\xi_{\underline{i},\underline{i}-nt_{\delta}(1,\cdots,1)}\xi_{\underline{i},\underline{j}+n(\varepsilon-e(\Omega)+e(\Omega)\delta-t_{\delta}(1,\cdots,1))}
\]
where the sum is over a representative set of all nontrivial
double cosets
$\Sigma_{\underline{i},e(\Omega)}\delta\Sigma_{\underline{i},\varepsilon-e(\Omega),\underline{j}}$
of $\Sigma_{\underline{i}}$,
$\mu_{\delta}=[\Sigma_{\underline{i},\underline{j},\varepsilon-e(\Omega)+e(\Omega)\delta}:\Sigma_{\underline{i},\underline{j},\varepsilon-e(\Omega),e(\Omega)\delta}]$
and $t_{\delta}$ is the minimal entry of
$\varepsilon-e(\Omega)+e(\Omega)\delta$. Now
$\varepsilon-e(\Omega)$ is of less volume. If $t_{\delta}=0$ then
$\varepsilon-e(\Omega)+e(\Omega)\delta$ is of the assuming form
and smaller than $\varepsilon$, and if $t_{\delta}>0$ then
$\varepsilon-e(\Omega)+e(\Omega)\delta-t_{\delta}(1,\cdots,1)$ is
of less volume. $\square$
\end{pf}

As a consequence we have
\begin{Thm}\label{T:noetherian} {\rm (i)} $\widetilde{S}(\underline{i})$
is a Noetherian ring.

{\rm (ii)} $\widetilde{S}(n,r)$ is a Noetherian ring.
\end{Thm}
\begin{pf} By Proposition~\ref{P:fg} both $\widetilde{S}(\underline{i})$ and $\widetilde{S}(n,r)$ are finitely generated over $\oplus_{\underline{i}\in I(n,r)/\Sigma_{r}}B_{\underline{i}}$, which is a Noetherian ring by Proposition~\ref{P:commsubalg}. $\square$
\end{pf}

In fact by a more subtle discussion we can reduce the number of
generators in Proposition~\ref{P:fg}. We assume
$\underline{i}=\underline{i}(1,\lambda_{1})\cdots\underline{i}(n,\lambda_{n})$.
Then

\begin{Cor} As a $B_{\underline{i}}$-module,
$K\{\xi_{\underline{i},\underline{j}+n\varepsilon}\ |\
\varepsilon\in\mathbb{Z}^{r}\}$ is generated by
\[\{\xi_{\underline{i},\underline{j}+n\varepsilon}\ |\ \theta^{0}_{s}(\varepsilon) {\text\rm\ is\ absolutely\ successive\ } \forall s=1,\cdots,n\}\]
\end{Cor}
\begin{pf} Denote by $M$ the module under investigation. Then
$M=M_{1}\boxtimes\cdots \boxtimes M_{n}$, where
$M_{s}=K\{\xi_{\underline{i},\underline{j}+n\varepsilon}\ |\
\varepsilon\in\Phi_{s}\}$ is a $B_{s}$-module. By
Proposition~\ref{P:fg}, $M_{s}$ is generated by
$\{\xi_{\underline{i},\underline{j}+n\varepsilon}\ |\
\varepsilon\in\Phi_{s},\theta^{0}_{s}(\varepsilon)$ is absolutely
successive $\forall s=1,\cdots,n\}$. $\square$
\end{pf}

\subsection{Special case : $\widetilde{S}(1\cdots r)$ ($n\geq r$)}

Assume $n\geq r$. Let $\underline{u}=(1,\cdots,r)\in I(n,r)$. Then
$\widetilde{S}(\underline{u})=K\{\xi_{\underline{u},\underline{u}w}\
|\ w\in\widehat{\Sigma}_{r}\}$ with
$\xi_{\underline{u},\underline{u}w}\xi_{\underline{u},\underline{u}w'}=\xi_{\underline{u},\underline{u}w'w}$
. Note that
$\xi_{\underline{u},\underline{u}w}=\xi_{\underline{u},\underline{u}w'}$
if and only if $w=w'$. Therefore the set of basis elements
$\{\xi_{\underline{u},\underline{u}w}|w\in\widehat{\Sigma}_{r}\}$
is closed under multiplication. In fact it is isomorphic to
$\widehat{\Sigma}_{r}$. Hence $\widetilde{S}(\underline{u})$ is
isomorphic to the group algebra $K\widehat{\Sigma}_{r}$. We will
identify these two algebras. Define {\em the Schur functor} $F$ :
$\widetilde{S}(n,r)$-$Mod$$\longrightarrow
K\widehat{\Sigma}_{r}$-$Mod$, $W\mapsto
\xi_{\underline{u},\underline{u}}W$.

\begin{Thm}\label{T:schurfunctor}
Assume $charK>r$ or $charK=0$.  Then the Schur functor $F$ is an
equivalence.
\end{Thm}
\begin{pf}
Since $\Sigma_{\underline{i}}\geq\Sigma_{\underline{u}}$ for any
$\underline{i}\in I(n,r)$, it follows by Lemma~\ref{L:summand}
that each $\widetilde{S}(n,r)\xi_{\underline{i},\underline{i}}$ is
a direct summand of
$\widetilde{S}(n,r)\xi_{\underline{u},\underline{u}}$. Therefore
$\widetilde{S}(n,r)\xi_{\underline{u},\underline{u}}$ is a
progenerator, and hence we have the desired equivalence. $\square$
\end{pf}

\subsection{Special case : $\widetilde{S}(112\cdots n)$ ($r=n+1$)}

Assume $r=n+1$. Let $\underline{v}=112\cdots n$.

\begin{Thm}
Assume $charK>n+1$ or $charK=0$. Then $\widetilde{S}(n,n+1)$ is
Morita equivalent to $\widetilde{S}(\underline{v})$.
\end{Thm}
\begin{pf}
For any $\underline{i}\in I(n,n+1)$ there exists $\sigma\in
\Sigma_{n+1}$ such that
$\Sigma_{\underline{i}\sigma}\geq\Sigma_{\underline{v}}$.
Therefore by Lemma~\ref{L:summand} each
$\widetilde{S}(n,n+1)\xi_{\underline{i},\underline{i}}$ is a
direct summand of
$\widetilde{S}(n,n+1)\xi_{\underline{v},\underline{v}}$. Therefore
$\widetilde{S}(n,n+1)\xi_{\underline{v},\underline{v}}$ is a
progenerator, and hence we have the desired Morita equivalence.
$\square$
\end{pf}

\section{The center}\label{S:center} In this section we will study the precise structure of $Z=Z(\widetilde{S}(n,r))$, the center
of $\widetilde{S}(n,r)$.

For $\varepsilon\in \mathbb{Z}^{r}$, set
\[ c_{\varepsilon}=\sum_{\underline{i}\in I(n,r)/\Sigma_{r}}\sum_{\sigma\in\Sigma_{\varepsilon}\backslash\Sigma_{r}/\Sigma_{\underline{i}}}\xi_{\underline{i},\underline{i}+n\varepsilon\sigma}\]
Then $c_{\varepsilon}=c_{\varepsilon'}$  if and only if
$\varepsilon\sim_{\Sigma_{r}}\varepsilon'$.
\begin{Prop}\label{P:center1} For any $\varepsilon\in\mathbb{Z}^{r}$,
$c_{\varepsilon}\in Z$.
\end{Prop}
\begin{pf}
Let $\underline{i},\underline{j}\in I(n,r)$,
$\varepsilon^{0}\in\mathbb{Z}^{r}$. Then by
Proposition~\ref{P:productformula2}
\[
\begin{array}{c}
\xi_{\underline{i},\underline{j}+n\varepsilon^{0}}c_{\varepsilon}=\xi_{\underline{i},\underline{j}+n\varepsilon^{0}}(\sum_{\sigma\in\Sigma_{\varepsilon}\backslash\Sigma_{r}/\Sigma_{\underline{j}}}\xi_{\underline{j},\underline{j}+n\varepsilon\sigma})
=\sum_{\sigma\in\Sigma_{\varepsilon}\backslash\Sigma_{r}/\Sigma_{\underline{j}}}\xi_{\underline{i},\underline{j}+n\varepsilon^{0}}\xi_{\underline{j},\underline{j}+n\varepsilon\sigma}\\
=\sum_{\sigma\in\Sigma_{\varepsilon}\backslash\Sigma_{r}/\Sigma_{\underline{j}}}\sum_{\delta\in\Sigma_{\underline{j},\varepsilon\sigma}\backslash\Sigma_{\underline{j}}/\Sigma_{\underline{i},\underline{j},\varepsilon^{0}}}[\Sigma_{\underline{i},\underline{j},\varepsilon\sigma\delta+\varepsilon^{0}}:\Sigma_{\underline{i},\underline{j},\varepsilon^{0},\varepsilon\sigma\delta}]\xi_{\underline{i},\underline{j}+n(\varepsilon\sigma\delta+\varepsilon^{0})}\\
=\sum_{w\in\Sigma_{\varepsilon}\backslash\Sigma_{r}/\Sigma_{\underline{i},\underline{j},\varepsilon^{0}}}[\Sigma_{\underline{i},\underline{j},\varepsilon w+\varepsilon^{0}}:\Sigma_{\underline{i},\underline{j},\varepsilon^{0},\varepsilon w}]\xi_{\underline{i},\underline{j}+n(\varepsilon w+\varepsilon^{0})}\\
\end{array}
\]
\[
\begin{array}{c}
c_{\varepsilon}\xi_{\underline{i},\underline{j}+n\varepsilon^{0}}=(\sum_{\sigma\in\Sigma_{\varepsilon}\backslash\Sigma_{r}/\Sigma_{\underline{i}}}\xi_{\underline{i},\underline{i}+n\varepsilon\sigma})\xi_{\underline{i},\underline{j}+n\varepsilon^{0}}
=\sum_{\sigma\in\Sigma_{\varepsilon}\backslash\Sigma_{r}/\Sigma_{\underline{i}}}\xi_{\underline{i},\underline{i}+n\varepsilon\sigma}\xi_{\underline{i},\underline{j}+n\varepsilon^{0}}\\
=\sum_{\sigma\in\Sigma_{\varepsilon}\backslash\Sigma_{r}/\Sigma_{\underline{i}}}\sum_{\delta\in\Sigma_{\underline{i},\varepsilon\sigma}\backslash\Sigma_{\underline{i}}/\Sigma_{\underline{i},\underline{j},\varepsilon^{0}}}[\Sigma_{\underline{i},\underline{j}, \varepsilon^{0}+\varepsilon\sigma\delta}:\Sigma_{\underline{i},\underline{j},\varepsilon\sigma\delta,\varepsilon^{0}}]\xi_{\underline{i},\underline{j}+n(\varepsilon^{0}+\varepsilon\sigma\delta)}\\
=\sum_{w\in\Sigma_{\varepsilon}\backslash\Sigma_{r}/\Sigma_{\underline{i},\underline{j},\varepsilon^{0}}}[\Sigma_{\underline{i},\underline{j},\varepsilon w+\varepsilon^{0}}:\Sigma_{\underline{i},\underline{j},\varepsilon^{0},\varepsilon w}]\xi_{\underline{i},\underline{j}+n(\varepsilon w+\varepsilon^{0})}\\
\end{array}\]

The last equalities in these two formulas follow from the
following lemma. $\square$
\end{pf}

\begin{Lem} Let $G$ be a group, $H_{1}$, $H_{2}$ are two subgroups
of $G$ and $H_{3}$ a subgroup of $H_{2}$. Let $\Sigma$ be a
representative set of $H_{1}\backslash G/H_{2}$, and for
$\sigma\in\Sigma$ let $\eta(\sigma)$ be a representative set of
$(H_{1}^{\sigma}\cap H_{2})\backslash H_{2}/H_{3}$, where
$H_{1}^{\sigma}=\sigma^{-1}H_{1}\sigma$. Then
$\cup_{\sigma\in\Sigma}\sigma\eta(\sigma)$ is a representative set
of $H_{1}\backslash G/H_{3}$.
\end{Lem}

\begin{Prop}\label{P:center2} $Z$ is spanned by $\{c_{\varepsilon}\ |\ \varepsilon\in\mathbb{Z}^{r}\}$.
\end{Prop}
Before proving this proposition, we first have a look at what form
a central element should have. Let $c=\sum_{\underline{i}\in
I(n,r)/\Sigma_{r}}\sum_{\underline{j}\in
I(n,r)/\Sigma_{\underline{i}}}\sum_{\varepsilon\in
\mathbb{Z}^{r}/\Sigma_{\underline{i},\underline{j}}}\lambda_{\underline{i},\underline{j},\varepsilon}\xi_{\underline{i},\underline{j}+n\varepsilon}
\in \widetilde{S}(n,r)$ be a central element. Then for any
$\underline{i}\in I(n,r)/\Sigma_{r}$, we have
$\xi_{\underline{i},\underline{i}}c\in
\widetilde{S}(\underline{i})$. Moreover
$\xi_{\underline{i},\underline{i}}c$ lies in the center of
$\widetilde{S}(\underline{i})$, and hence lies in
$B_{\underline{i}}$ by Proposition~\ref{P:bimaximal}. Therefore
$\lambda_{\underline{i},\underline{j},\varepsilon}=0$ if
$\underline{j}\neq \underline{i}$. Namely $c$ can be written as
\[c=\sum_{\underline{i}\in
I(n,r)/\Sigma_{r}}\ \sum_{\varepsilon\in
\mathbb{Z}^{r}/\Sigma_{\underline{i}}}\lambda_{\underline{i},\varepsilon}\xi_{\underline{i},\underline{i}+n\varepsilon}
\]

\begin{Lem}\label{L:bridge} Let $\underline{i}\in I(n,r)$ and assume $\xi\in B_{\underline{i}}$ satisfies
$\xi_{1\cdots 1,\underline{i}}\xi=0$. Then $\xi=0$.
\end{Lem}
\begin{pf} Suppse
$\xi=\sum_{\varepsilon\in\mathbb{Z}^{r}/\Sigma_{\underline{i}}}\lambda_{\varepsilon}\xi_{\underline{i},\underline{i}+n\varepsilon}\in
B_{\underline{i}}$. Then $\xi_{1\cdots
1,\underline{i}}\xi=\sum_{\varepsilon\in\mathbb{Z}^{r}/\Sigma_{\underline{i}}}\lambda_{\varepsilon}\xi_{1\cdots
1,\underline{i}+n\varepsilon}$. But for
$\varepsilon,\varepsilon'\in\mathbb{Z}^{r}/\Sigma_{\underline{i}}$,
we have $\xi_{1\cdots 1,\underline{i}+n\varepsilon}=\xi_{1\cdots
1,\underline{i}+n\varepsilon'}$ if and only if
$\varepsilon=\varepsilon'$. Therefore $\xi_{1\cdots
1,\underline{i}}\xi=0$ implies $\lambda_{\varepsilon}=0$ for any
$\varepsilon\in\mathbb{Z}^{r}/\Sigma_{\underline{i}}$. $\square$
\end{pf}
Now we are ready for the\\
{\bf Proof of Proposition~\ref{P:center2}\ : }

Let $c=\sum_{\underline{i}\in
I(n,r)/\Sigma_{r}}\sum_{\varepsilon\in
\mathbb{Z}^{r}/\Sigma_{\underline{i}}}\lambda_{\underline{i},\varepsilon}\xi_{\underline{i},\underline{i}+n\varepsilon}\in\widetilde{S}(n,r)$
be a central element. Then
$\xi=c-\sum_{\varepsilon\in\mathbb{Z}^{r}/\Sigma_{\underline{i}}}\lambda_{1\cdots
1,\varepsilon}c_{\varepsilon}$ is also a central element. Moreover
$\xi_{1\cdots
1,\underline{i}}\xi\xi_{\underline{i},\underline{i}}=\xi_{1\cdots
1,\underline{i}}\xi=\xi\xi_{1\cdots 1,\underline{i}}=0$ for any
$\underline{i}\in I(n,r)$. It follows from Lemma~\ref{L:bridge}
that $\xi\xi_{\underline{i},\underline{i}}=0$ for any
$\underline{i}\in I(n,r)$, and hence $\xi=0$. Therefore
$c=\sum_{\varepsilon\in\mathbb{Z}^{r}/\Sigma_{\underline{i}}}\lambda_{1\cdots
1,\varepsilon}c_{\varepsilon}$, as desired. $\square$

\begin{Thm}\label{T:center} The center $Z$ of $\widetilde{S}(n,r)$ is isomorphic to $K[t_{1},\cdots,t_{r-1},t_{r},t_{r}^{-1}]$, where $t_{1},\cdots,t_{r}$ are indeterminates. To be precise,
$Z=K[c_{\varepsilon^{1}},\cdots,c_{\varepsilon^{r-1}},c_{\varepsilon^{r}},c_{-\varepsilon^{r}}]$,
where $\varepsilon^{1},\cdots,\varepsilon^{r}$ are defined in
Section~\ref{S:1...1}. In particular, $\widetilde{S}(n,r)$ is
indecomposable.
\end{Thm}
\begin{pf} Sending
$c_{\varepsilon}$ to $c_{\varepsilon}\xi_{1\cdots 1,1\cdots 1}$
defines a $K$-algebra isomorphism from $Z$ to
$\widetilde{S}(1\cdots 1)$. $\square$
\end{pf}

\section{Examples}\label{S:examples}

We shall denote by $\mathcal{M}_{n}$ the algebra of $n\times
n$-matrices with entries from the field $K$.

\begin{eg}
Let $r=1$. Then $\widetilde{S}(n,1)$ is isomorphic to
$\mathcal{M}_{n}\otimes K[t,t^{-1}]$. So $\widetilde{S}(n,1)$ is
Morita equivalent to $K[t,t^{-1}]= K\widehat{\Sigma}_{1}$ with the
equivalence given by the Schur functor. The center of this algebra
is $1\otimes K[t,t^{-1}]$.
\end{eg}

\begin{eg}
Let $n=1$. Then by Proposition~\ref{P:1...1}
$\widetilde{S}(1,r)\cong
K[t_{1},\cdots,t_{r-1},t_{r},t_{r}^{-1}]$, where
$t_{1},\cdots,t_{r}$ are indeterminates. This is a commutative
algebra.
\end{eg}

\begin{eg} Let
$n=2$, $r=2$. Then by Theorem~\ref{T:center} the center of
$\widetilde{S}(2,2)$ is $Z=K[c_{10},c_{11},c_{11}^{-1}]$. We shall
classify all simple $\widetilde{S}(2,2)$-modules. Let $M$ be a
simple $\widetilde{S}(2,2)$-module, then $Z$ acts as scalars, say,
$c_{10}$ as $a\in K$, and $c_{11}$ as $b\in K^{\times}$. Then the
action of $\widetilde{S}(2,2)$ on $M$ factors through the algebra
$A_{a,b}=\widetilde{S}(2,2)/(c_{10}-a,c_{11}-b)$.

Case 1 : $charK\neq 2$. As a $Z$-module, $\widetilde{S}(2,2)$ is
free with basis
\[\tiny
\begin{array}{cccc}
\xi_{11,11} & \xi_{11,12} & \xi_{11,12}\gamma & \xi_{11,22}\\
\frac{\xi_{12,11}}{2} & \frac{\xi_{12,12}+\xi_{12,21}}{2} & \gamma
&
\frac{\xi_{12,22}}{2}\\
\frac{\gamma'\xi_{12,11}}{2} & \gamma' &
\frac{\xi_{12,12}-\xi_{12,21}}{2} &
\frac{\gamma'\xi_{12,22}}{2}\\
\xi_{22,11} & \xi_{22,12} & \xi_{22,12}\gamma & \xi_{22,22}
\end{array}
\]
where $\frac{\xi_{12,12}-\xi_{12,21}}{2}\ \gamma'\
\frac{\xi_{12,12}+\xi_{12,21}}{2}=\gamma'$,
$\frac{\xi_{12,12}+\xi_{12,21}}{2}\ \gamma\
\frac{\xi_{12,12}-\xi_{12,21}}{2}=\gamma$,
$\gamma'\gamma=(c_{10}^{2}-4c_{11})\frac{\xi_{12,12}-\xi_{12,21}}{2}$
and
$\gamma\gamma'=(c_{10}^{2}-4c_{11})\frac{\xi_{12,12}+\xi_{12,21}}{2}$.
Other structure constants are easy to calculate.

When $a^{2}-4b=0$, the radical $radA_{a,b}$ of $A_{a,b}$ is
spanned by $\xi_{11,12}\gamma$, $\gamma$,
$\frac{\gamma'\xi_{12,11}}{2}$, $\gamma'$,
$\frac{\gamma'\xi_{12,22}}{2}$, and $\xi_{22,12}\gamma$, and hence
 $A_{a,b}/rad(A_{a,b})$ is isomorphic to
$\mathcal{M}_{3}\times K$. When $a^{2}-4b\not=0$, the algebra
$A_{a,b}$ is isomorphic to $\mathcal{M}_{4}$.

Therefore the isoclasses of simple $\widetilde{S}(2,2)$-modules
are parametrized by the plane without the $x$-axis and with the
curve $\{(x,\frac{x^{2}}{4})\ |\ x\neq 0\}$ doubled.

Case 2 : $charK=2$. As a $Z$-module, $\widetilde{S}(2,2)$ is free
with basis
\[\tiny
\begin{array}{cccc}
\xi_{11,11} & \xi_{11,12} & \xi_{11,12}l & \xi_{11,22}\\
\xi_{12,11} & \xi_{12,12} & \xi_{12,21} & \xi_{12,22}\\
{ l\xi_{12,11}} & l & l\xi_{12,21} & l\xi_{12,22}\\
 \xi_{22,11} & \xi_{22,12} & \xi_{22,12}l &
\xi_{22,22}
\end{array}
\]
where $\xi_{12,12}l=l\xi_{12,12}=l$,
$l^{2}=c_{10}l+c_{11}\xi_{12,12}$,
$\xi_{12,21}l=l\xi_{12,21}+c_{10}\xi_{12,21}$,
$\xi_{11,12}l\xi_{12,11}=c_{10}\xi_{11,11}$,
$\xi_{11,12}l\xi_{12,22}=c_{10}\xi_{11,22}$,
$\xi_{22,12}l\xi_{12,11}=c_{10}\xi_{22,11}$, and
$\xi_{22,12}l\xi_{22,22}=c_{10}\xi_{12,22}$. Other structure
constants are easy to calculate.

We denote the images of these elements in $A_{a,b}$ by the same
notations. For simplicity we assume $K$ is quadratically closed.
Let $\alpha=\xi_{12,12}+\xi_{12,21}$,
$\beta=l\xi_{12,21}+\sqrt{b}\xi_{12,12}$, then $A_{a,b}$ has the
following basis
\[\tiny\begin{array}{cccc}
\xi_{11,11} & \xi_{11,12} & \xi_{11,12}\beta & \xi_{11,22}\\
\xi_{12,11} & \xi_{12,12} & \alpha & \xi_{12,22}\\
\beta\xi_{12,11} & \beta & \alpha\beta & \beta\xi_{12,22}\\
\xi_{22,11} & \xi_{22,12} & \xi_{22,12}\beta &
\xi_{22,22}\end{array}\] where $\beta^{2}=0$,
$\beta\alpha=\alpha\beta+a$, $\xi_{12,12}\beta\xi_{12,12}=\beta$,
$\xi_{\underline{i},12}\beta\xi_{12,\underline{j}}=a\xi_{\underline{i},\underline{j}}$,
$\xi_{12,\underline{i}}\xi_{\underline{i},12}=\alpha$,
$\xi_{\underline{i},12}\xi_{12,\underline{j}}=0$. In the above,
$\underline{i},\underline{j}\in \{11,22\}$.

If $a=0$, then the radical $radA_{0,b}$ of $A_{0,b}$ is spanned by
$\xi_{11,12}$, $\xi_{11,12}\beta$, $\xi_{12,11}$, $\alpha$,
$\xi_{12,22}$, $\beta\xi_{12,11}$, $\beta$, $\alpha\beta$,
$\beta\xi_{12,22}$, $\xi_{22,12}$, $\xi_{22,12}\beta$. Therefore
$A_{0,b}/radA_{0,b}$ is isomorphic to $\mathcal{M}_{2}\times K$.
If $a\neq 0$, then the algebra $A_{a,b}$ has basis
\[\tiny\begin{array}{cccc}
\xi_{11,11} & \xi_{11,12} & \xi_{11,12}\beta & \xi_{11,22}\\
a^{-1}\beta\xi_{12,11} & 1+a^{-1}\alpha\beta & \beta & a^{-1}\beta\xi_{12,22}\\
a^{-1}\xi_{12,11} & a^{-1}\alpha & a^{-1}\alpha\beta & a^{-1}\xi_{12,22}\\
\xi_{22,11} & \xi_{22,12} & \xi_{22,12}\beta &
\xi_{22,22}\end{array}\] If we denote by $E_{ij}$ the element in
the $(i,j)$-entry then $E_{ij}E_{kl}=\delta_{jk}E_{il}$. In
particular $A_{a,b}$ is isomorphic to $\mathcal{M}_{4}$.

Therefore the isoclasses of simple $\widetilde{S}(2,2)$-modules
are parametrized by the plane without $x$-axis and with the curve
$\{(0,y)\ |\ y\neq 0\}$ doubled.

\end{eg}

\bibliographystyle{amsplain}

{Dong Yang}

{Department of Mathematical Sciences, Tsinghua University,
Beijing100084, P.R.China.}

 {\it Email address : }
{yangdong98@mails.tsinghua.edu.cn.}

\end{document}